\documentclass[a4paper,12pt]{article}
\usepackage{authblk}
\usepackage{amsmath}
\usepackage{amsfonts}
\usepackage{amssymb}
\usepackage{amsthm}
\usepackage{graphicx}
\theoremstyle{remark}

\begin{document}
\title{On the characterization of drilling rotation in the 6--parameter resultant shell theory}
\author{Mircea B\^{\i}rsan%
\thanks{\,Mircea B\^{\i}rsan,  Fakult\"at f\"ur  Mathematik, Universit\"at Duisburg-Essen, Campus Essen, Thea-Leymann Str. 9, 45127 Essen, Germany, email: mircea.birsan@uni-due.de ; and
Department of Mathematics, University ``A.I. Cuza'' of Ia\c{s}i, 700506 Ia\c{s}i,  Romania}\,\,\,
and\, Patrizio Neff\,%
\thanks{Patrizio Neff, Lehrstuhl f\"ur Nichtlineare Analysis und Modellierung, Fakult\"at f\"ur  Mathematik, Universit\"at Duisburg-Essen, Campus Essen, Thea-Leymann Str. 9, 45127 Essen, Germany, email: patrizio.neff@uni-due.de, Tel.: +49-201-183-4243}
}


\maketitle

\begin{abstract}
We analyze  geometrically non-linear isotropic elastic shells and prove the existence of minimizers. In general, the model takes into account the effect of drilling rotations in shells. For the special case of shells without drilling rotations we present a representation theorem for the strain energy function.
\end{abstract}

\section{Introduction}\label{sec1}
The paper is concerned with the geometrically non-linear 6-parameter resultant shell theory. This model of shells involves two independent kinematic fields: the translation vector field and the rotation tensor field, which have in total 6 independent scalar kinematic variables. This shell theory   was originally proposed by Reissner \cite{Reissner74} and was   developed consistently by several authors \cite{Libai98,Pietraszkiewicz-book04}.

In Section \ref{sec2} we briefly present the kinematics of 6-parameter shells, as well as the equations of equilibrium and the constitutive equations of elastic shells. We formulate the boundary-value problem and prove the existence of minimizers associated to the deformation of isotropic shells in Section \ref{sec3}. This model is able to describe the effect of drilling rotations in shells. In order to see the difference to the Reissner-type kinematics for shells, we analyze in Section \ref{sec4} shells without drilling rotations and give a representation theorem for the strain energy function. Finally, we consider isotropic shells without drilling rotations and identify the constitutive coefficients by comparison with the 6-parameter model. In case of shells without drilling rotations the strain energy function is only positive semi-definite, so that the general theorem for the existence of minimizers does not apply.

\section{Equations of equilibrium}\label{sec2}
\vspace{1pt}

We consider a shell with the base surface $S^0$ in the reference cofiguration characterized by the position vector (relative to a fixed point $O$) $\boldsymbol{y}^0:\omega\subset \mathbb{R}^2\rightarrow\mathbb{R}^3$, $\boldsymbol{y}^0=\boldsymbol{y}^0(x_1,x_2)$, and the structure tensor $\boldsymbol{Q}^{0}:\omega\subset \mathbb{R}^2\rightarrow SO(3)$, $ \boldsymbol{Q}^{0} = \boldsymbol{d}_i^0 (x_1,x_2)\otimes \boldsymbol{e}_i\,$. Here, $(x_1,x_2)$ are the material curvilinear coordinates on $S^0\,$, $\{\boldsymbol{e}_i\}$ is the fixed orthonormal vector basis of the Euclidean space, and $\{\boldsymbol{d}^0_i(x_1,x_2)\}$ is the orthonormal triad of directors which characterizes the orthogonal tensor field $\boldsymbol{Q}^{0} $ \cite{Pietraszkiewicz-book04,Eremeyev06}. We employ the usual notations: the Latin indexes $i,j,...$ take the values $\{1,2,3\}$, the Greek indexes $\alpha,\beta,...$ the values $\{1,2\}$, the partial derivative $\partial_\alpha f= {\partial f}/{\partial x_\alpha}\,$, as well as the Einstein summation convention over repeated indexes. The set $\omega$  is assumed to be a bounded open domain  with Lipschitz boundary in the $Ox_1x_2$ plane.

For the deformed configuration of the shell, we denote by $S$ the base surface, $\boldsymbol{y}(x_1,x_2)$ the position vector and $\{\boldsymbol{d}_i(x_1,x_2)\}$ the orthonormal triad of directors. The displacement vector is defined as usual by $\boldsymbol{u}= \boldsymbol{y}-\boldsymbol{y}^0\,$ and the elastic rotation (between $S^0$ and $S$) by the  proper orthogonal tensor field $\boldsymbol{Q}= \boldsymbol{d}_i\otimes \boldsymbol{d}_i^0\,$. The orthogonal tensor $\boldsymbol{R}=\boldsymbol{Q}\boldsymbol{Q}^{0}= \boldsymbol{d}_i\otimes \boldsymbol{e}_i$ describes the total rotation from $\omega$ to $S$.

Let $\boldsymbol{a}_\alpha= \partial_\alpha \boldsymbol{y}^0$ be the (covariant) base vectors in the tangent plane to $S^0$,  $\boldsymbol{n}^0=\boldsymbol{a}_1\times\boldsymbol{a}_2 /{ \| \boldsymbol{a}_1\times\boldsymbol{a}_2\|}$ the unit normal to $S^0$, and $\{\boldsymbol{a}^\alpha\}$ the reciprocal (contravariant) basis in the tangent plane, with $\boldsymbol{a}_\alpha\cdot \boldsymbol{a}_\beta=a_{\alpha\beta}\,$ and $ \boldsymbol{a}^\alpha\cdot\boldsymbol{a}_\beta =\delta^\alpha_\beta$ (the Kronecker symbol).
Then the shell deformation gradient tensor is expressed by  $\boldsymbol{F}=\text{Grad}_s\boldsymbol{y}=\partial_\alpha \boldsymbol{y}\otimes\boldsymbol{a}^\alpha$.

We designate by $\boldsymbol{N}$ and $\boldsymbol{M}$   the internal surface stress resultant and stress couple tensors of the 1$^{st}$ Piola--Kirchhoff type for the shell, and by $\boldsymbol{f}$ and $\boldsymbol{c}$   the external surface resultant force and couple vectors applied to points of $S$, but measured per unit area of $S^0\,$.
The equilibrium equations for 6-parameter shells are \cite{Eremeyev06}
\begin{equation}\label{1}
\begin{array}{l}
    \mathrm{Div}_s\, \boldsymbol{N}+\boldsymbol{f}=\boldsymbol{0},\\
    \mathrm{Div}_s\, \boldsymbol{M} + \mathrm{axl}(\boldsymbol{N}\boldsymbol{F}^T-\boldsymbol{F}\boldsymbol{N}^T)
    +\boldsymbol{c}=\boldsymbol{0},
    \end{array}
\end{equation}
where Div$_s$ is the surface divergence,  $(\cdot)^T$ denotes the transpose, and axl$(\,\cdot)$ represents the axial vector of a skew--symmetric tensor. We consider boundary conditions of the type \cite{Pietraszkiewicz11}
\begin{equation}\label{2}
\begin{array}{l}
\boldsymbol{N}\boldsymbol{\nu}=\boldsymbol{n}^*,\qquad \boldsymbol{M}\boldsymbol{\nu}=\boldsymbol{m}^*\qquad\mathrm{along}\,\,\,\partial S^0_f\,,\\
    \quad\boldsymbol{y}=\boldsymbol{y}^* ,\qquad\quad\,\, \boldsymbol{R}=\boldsymbol{R}^* \qquad\mathrm{along}\,\,\,\partial S^0_d\,,
    \end{array}
\end{equation}
where $\boldsymbol{\nu}$ is the external unit normal vector to the boundary curve $\partial S^0$ (lying in the tangent plane) and $\{ \partial S^0_f\,, \partial S^0_d\,\}$ is a disjoint partition of $\partial S^0\,$.

According to \cite{Eremeyev06,Pietraszkiewicz-book04}, the elastic shell strain tensor $\boldsymbol{E}^e$ and the
bending--curvature  tensor $\boldsymbol{K}^e$
in the material representation are
\begin{equation}\label{3}
\begin{array}{l}
    \boldsymbol{E}^e=\boldsymbol{Q}^{ T}\!\boldsymbol{F}-\! \text{Grad}_s\,\boldsymbol{y}^0=\big(\boldsymbol{Q}^{ T}\partial_\alpha \boldsymbol{y}- \!\boldsymbol{a}_\alpha \big) \!\otimes \boldsymbol{a}^\alpha\\
    \boldsymbol{K}^e=  \text{axl}(\boldsymbol{Q}^{T} \partial_\alpha \boldsymbol{Q} ) \otimes \boldsymbol{a}^\alpha .
\end{array}
\end{equation}
Under the hyperelasticity assumption, $\boldsymbol{N}$ and $\boldsymbol{M}$ are expressed by the constitutive equations
\begin{equation}\label{4}
    \boldsymbol{N}=\boldsymbol{Q}\,\dfrac{\partial\, W}{\partial \boldsymbol{E}^e}\,\,,\qquad \boldsymbol{M}=\boldsymbol{Q}\,\dfrac{\partial\, W}{\partial \boldsymbol{K}^e}\,\,,
\end{equation}
where $W=W(\boldsymbol{E}^e,\boldsymbol{K}^e)$ is the strain energy density of the elastic shell.
The boundary--value problem describing the deformation of non-linear elastic shells consists of the equations \eqref{1}-\eqref{4}. We assume the existence of a function $\Lambda(\boldsymbol{y},\boldsymbol{R})$ representing the potential of external surface loads $\boldsymbol{f}$, $\boldsymbol{c}$, and boundary loads $\boldsymbol{n}^*$, $\boldsymbol{m}^*$  \cite{Pietraszkiewicz04}. This loading potential can be decomposed additively as
\begin{equation*}
\begin{array}{l}
    \Lambda(\boldsymbol{y},\boldsymbol{R})=\Lambda_{S^0}(\boldsymbol{y},\boldsymbol{R}) + \Lambda_{\partial S^0_f}(\boldsymbol{y},\boldsymbol{R}),\\
    \Lambda_{S^0}(\boldsymbol{y},\boldsymbol{R})= \displaystyle{\int_{S^0}}\! \boldsymbol{f}\!\cdot\! \boldsymbol{u}\, \mathrm{d}S + \Pi_{S^0}(\boldsymbol{R}), \\
    \Lambda_{\partial S^0_f}(\boldsymbol{y},\boldsymbol{R})= \displaystyle{\int_{\partial S^0_f}}\! \boldsymbol{n}^*\!\cdot\! \boldsymbol{u}\, \mathrm{d}l + \Pi_{\partial S^0_f}(\boldsymbol{R}),
    \end{array}
\end{equation*}
where the load potential functions $\,\,\Pi_{S^0}\,,  \,\Pi_{\partial S^0_f}:\boldsymbol{L}^2( \omega ,SO(3))\rightarrow\mathbb{R}$ are   assumed to be continuous and bounded operators.
Corresponding to the deformation of elastic shells, we consider the following two--field minimization problem: find the pair $(\hat{\boldsymbol{y}},\hat{\boldsymbol{R}})$ in the admissible set $\mathcal{A}$ which realizes the minimum of the functional
\begin{equation}\label{5}
I(\boldsymbol{y},\boldsymbol{R})=\int_{S^0} W(\boldsymbol{E}^e,\boldsymbol{K}^e)\,\mathrm{d}S - \Lambda(\boldsymbol{y},\boldsymbol{R})
\end{equation}
for $(\boldsymbol{y},\boldsymbol{R})\in \mathcal{A}$, where $    \mathcal{A}:=\big\{(\boldsymbol{y},\boldsymbol{R})\in\boldsymbol{H}^1(\omega, \mathbb{R}^3)\times\boldsymbol{H}^1(\omega, SO(3))\,\big|\,\, \boldsymbol{y}_{| \partial S^0_d}=\boldsymbol{y}^*, \,\boldsymbol{R}_{| \partial S^0_d}=\boldsymbol{R}^* \big\}$.
Here, the boundary conditions are to be understood in the sense of traces, $\boldsymbol{H}^1$ denotes as usual the Sobolev space, and $\boldsymbol{L}^2$ the Lebesgue space. The variational principle associated to the total energy of elastic shells   \eqref{5} has been proved in \cite{Pietraszkiewicz04}.

\section{Shells with drilling rotations: existence of \\minimizers}
\label{sec3}

In case of physically linear isotropic shells, the strain energy density is assumed as the quadratic form
\begin{equation}\label{6}
    \begin{array}{l}
   2 W(\boldsymbol{E}^e,\boldsymbol{K}^e)=  \alpha_1\big(\mathrm{tr}  \boldsymbol{E}^e_{\parallel}\big)^2 +\alpha_2\,  \mathrm{tr} \big(\boldsymbol{E}^e_{\parallel}\big)^2    
   + \alpha_3\, \mathrm{tr}\big(\boldsymbol{E}^{e,T}_{\parallel}  \boldsymbol{E}^e_{\parallel} \big)  + \alpha_4     (\boldsymbol{n}^0 \boldsymbol{E}^e)^2\\
   \qquad \qquad \qquad  + \beta_1\big(\mathrm{tr}  \boldsymbol{K}^e_{\parallel}\big)^2 
       +\beta_2 \, \mathrm{tr} \big(\boldsymbol{K}^e_{\parallel}\big)^2    + \beta_3 \, \mathrm{tr}\big(\boldsymbol{K}^{e,T}_{\parallel} \boldsymbol{K}^e_{\parallel} \big)  + \beta_4     (\boldsymbol{n}^0 \boldsymbol{K}^e)^2,
\end{array}
\end{equation}
where $\boldsymbol{E}^e_{\parallel}\!=\boldsymbol{E}^e\!-\! \boldsymbol{n}^0\! \otimes\!\boldsymbol{n}^0 \!\boldsymbol{E}^e,\, \boldsymbol{K}^e_{\parallel}\!=\boldsymbol{K}^e\!-\! \boldsymbol{n}^0\!\otimes\!\boldsymbol{n}^0\! \boldsymbol{K}^e$.
The eight   coefficients $\alpha_k\,$, $\beta_k$   can depend in general on the initial structure curvature tensor  $ \boldsymbol{K}^0=\text{axl}\big(\partial_\alpha \boldsymbol{Q}^{0}\boldsymbol{Q}^{0,T}\big) \otimes \boldsymbol{a}^\alpha$. For the sake of simplicity, we assume in our discussion that $\alpha_k$ and $\beta_k$ are constant.\smallskip

\textbf{Theorem 1.} \textit{ Assume that the initial position vector $\boldsymbol{y}^0$ is continuous and injective and}
\begin{equation*}
    \begin{array}{l}
   \boldsymbol{y}^0\in\boldsymbol{H}^1(\omega ,\mathbb{R}^3),\qquad \boldsymbol{Q}^{0}\in\boldsymbol{H}^1(\omega, SO(3)),\\
   \partial_\alpha \boldsymbol{y}^0\in \boldsymbol{L}^\infty(\omega ,\mathbb{R}^3),\quad  \det\big(a_{\alpha\beta}(x_1,x_2)\big)\geq a_0^2 >0,
    \end{array}
\end{equation*}
\textit{where $a_0$ is a constant. The external loads and boundary data are assumed to satisfy the conditions}
\begin{equation*}
 \begin{array}{l}
    \boldsymbol{f}\in\boldsymbol{L}^2(\omega,\mathbb{R}^3),\qquad  \boldsymbol{n}^*\in \boldsymbol{L}^2(\partial\omega_f,\mathbb{R}^3),\\
    \boldsymbol{y}^*\in\boldsymbol{H}^1(\omega ,\mathbb{R}^3),\quad \, \boldsymbol{R}^*\in\boldsymbol{H}^1(\omega, SO(3)).
\end{array}
\end{equation*}
\textit{If the constitutive coefficients satisfy the conditions}
\begin{equation}\label{7}
    \begin{array}{l}
    2\alpha_1+\alpha_2+\alpha_3>0,\qquad \alpha_2+\alpha_3>0,\\
    \alpha_3-\alpha_2>0, \quad\, \alpha_4>0,\quad\,
    2\beta_1+\beta_2+\beta_3>0,\\
     \beta_2+\beta_3>0,\qquad \beta_3-\beta_2>0,\qquad  \beta_4>0,
\end{array}
\end{equation}
\textit{then the minimization problem \eqref{5} admits at least one minimizing solution pair
$(\hat{\boldsymbol{y}},\hat{\boldsymbol{R}})\in  \mathcal{A}$.}\smallskip

\textbf{Proof.} In view of the inequalities \eqref{7} we can deduce that there exists a constants $C_0>0$ such that
\begin{equation*}
    W( {\boldsymbol{E}^e}, {\boldsymbol{K}^e})\,\geq\, C_0\,\big( \,   \|\boldsymbol{E}^e\|^2 +  \|\boldsymbol{K}^e\|^2\,\big).
\end{equation*}
Moreover, in this case the strain energy density $W(\boldsymbol{E}^e,\boldsymbol{K}^e)$ is a strictly  convex function of its arguments. Then, according to   Theorem 6 from \cite{Birsan-Neff-MMS-2013}, we derive the existence of minimizers. The proof is based on the direct methods of the calculus of variations.$\hfill\Box$\smallskip

\textbf{Remark 2.} For isotropic  shells, the simplest expression of $W(\boldsymbol{E}^e,\boldsymbol{K}^e)$ corresponds to the form \eqref{6} with
\begin{equation}\label{8}
    \begin{array}{l}
    \alpha_1=C\nu,\qquad \alpha_2=0,\qquad \alpha_3=C(1-\nu),\\
     \alpha_4=\alpha_sC(1-\nu), \qquad     \beta_1=D\nu,\qquad \beta_2=0,\\
      \beta_3=D(1-\nu), \qquad \beta_4=\alpha_tD(1-\nu),
\end{array}
\end{equation}
where $h$ is the thickness of the shell, $E $ the Young modulus, $ \nu$ the  Poisson ratio of the material, $C= {Eh}/{(1-\nu^2)}$ is the stretching (in-plane) stiffness, $D= {E\,h^3}/{12(1-\nu^2)}$ is the bending stiffness, and $\alpha_s\,$, $\alpha_t$ are two shear correction factors. Note that the conditions \eqref{7} are fulfilled for the coefficients \eqref{8}. $\hfill\Box$\smallskip

Without loss of generality, one can choose the directors $\{\boldsymbol{d}_i^0\}$ such that $\boldsymbol{d}_3^0=\boldsymbol{n}^0$ is the unit normal to $S^0$. In what follows, we assume that $\boldsymbol{d}_3^0=\boldsymbol{n}^0$.

In the 6-parameter shell theory the drilling rotations are taken into account. The drilling rotation in a given point $S$ can be interpreted as the rotation about the director $\boldsymbol{d}_3\,$. The general form of rotations  about  $\boldsymbol{d}_3 $ is
\begin{equation*} 
    \begin{array}{l}
    \boldsymbol{R}_\theta=\boldsymbol{d}_3\otimes\boldsymbol{d}_3+ \cos \theta(\boldsymbol{1}\!-\!\boldsymbol{d}_3\otimes\boldsymbol{d}_3)+\sin\theta (\boldsymbol{d}_3\!\times\! \boldsymbol{1}),
\end{array}
\end{equation*}
where $\theta=\theta(x_1,x_2)$ is the rotation angle and $\boldsymbol 1= \boldsymbol{d}_i\otimes \boldsymbol{d}_i$ is the unit tensor.

Let us describe next  shells without drilling rotations in the framework of the 6-parameter shell theory.

\section{Shells without drilling rotations:\\ characterization}\label{sec4}

In case of shells without drilling rotations the strain energy density $W$ must be insensible  to the superposition of rotations $\boldsymbol{R}_\theta $ about $\boldsymbol{d}_3 \,$. This means that $ W( {\boldsymbol{E}^e}, {\boldsymbol{K}^e})$ is assumed to remain invariant under the transformation
\begin{equation}\label{10}
    \begin{array}{c}
      \boldsymbol{Q}\,\,\,\longrightarrow\,\,\, \boldsymbol{R}_\theta \boldsymbol{Q}\,\,.
\end{array}
\end{equation}
In view of the definitions \eqref{3}, this is equivalent to
\begin{equation}\label{11}
    \begin{array}{c}
      W( {\boldsymbol{E}^e}, {\boldsymbol{K}^e})=
      W\big([\boldsymbol{Q}^{ T}\boldsymbol{R}_\theta^T\partial_\alpha \boldsymbol{y}- \!\boldsymbol{a}_\alpha ] \!\otimes \boldsymbol{a}^\alpha, 
      \,\,\text{axl}[\boldsymbol{Q}^{T} \boldsymbol{R}_\theta^T \partial_\alpha (\boldsymbol{R}_\theta\boldsymbol{Q} )] \otimes \boldsymbol{a}^\alpha\big)
\end{array}
\end{equation}
for any angle $\theta(x_1,x_2)$. The following result gives a characterization of shells without drilling rotation.\smallskip

\textbf{Theorem 3.} \textit{ Assume that the strain energy function $W$ is invariant under the transformation \eqref{10}. Then, $W$ can be represented as a function of the arguments}
\begin{equation}\label{12}
    \begin{array}{r}
      W =
      \widetilde W\big(\boldsymbol{F}^{ T}\boldsymbol{F} \,\,, \,\,
       \boldsymbol{d}_3 \boldsymbol{F} \,\,, \,\,  \boldsymbol{F}^T\mathrm{Grad}_s\boldsymbol{d}_3 \big).
\end{array}
\end{equation}
\textit{Conversely, any  function $W$ of the form \eqref{12} is invariant under the superposition of drilling rotations \eqref{10}}.\smallskip

\textbf{Proof.} Firstly, it is clear that the function \eqref{12} is invariant under the drill rotation, since $ \boldsymbol{F}=\mathrm{Grad}_s\boldsymbol{y}$ and $\mathrm{Grad}_s\boldsymbol{d}_3$ are both independent of rotations about $ \boldsymbol{d}_3\,$. Conversely, let us assume that $W$ is invariant under the transformation \eqref{10}. If we denote by
$\boldsymbol{d}_1^\theta= \boldsymbol{R}_\theta\boldsymbol{d}_1=\cos\theta\boldsymbol{d}_1+ \sin\theta\boldsymbol{d}_2\,$ and $
    \boldsymbol{d}_2^\theta= \boldsymbol{R}_\theta\boldsymbol{d}_2=-\sin\theta\boldsymbol{d}_1+ \cos\theta\boldsymbol{d}_2\,$,
then we find $\boldsymbol{R}_\theta\boldsymbol{Q}=\boldsymbol{d}_i^\theta\otimes\boldsymbol{d}_i^0\,$. Inserting the last relation into equation \eqref{11} and imposing the conditions that the derivative of \eqref{11} with respect to $\theta$ and $\partial_\alpha\theta$ are zero, we obtain the equations
\begin{equation}\label{13}
    \begin{array}{l}
    \dfrac{\partial W}{\partial \boldsymbol{E}^e}\,\cdot\boldsymbol{c}(\boldsymbol{E}^e+\boldsymbol{a})+ \dfrac{\partial W}{\partial \boldsymbol{K}^e}\,\cdot\boldsymbol{c}(\boldsymbol{K}^e+\boldsymbol{K}^0)=0\\
    \mathrm{and}\qquad \dfrac{\partial W}{\partial(\boldsymbol{n}^0 \boldsymbol{K}^e)}\,=\boldsymbol{0},
\end{array}
\end{equation}
where we have used the notations $\boldsymbol{a}=\boldsymbol{a}_\alpha\otimes\boldsymbol{a}^\alpha= \boldsymbol{d}_\alpha^0\otimes\boldsymbol{d}^0_\alpha$ and $\boldsymbol{c}=  \boldsymbol{d}_1^0\otimes\boldsymbol{d}^0_2-  \boldsymbol{d}_2^0\otimes\boldsymbol{d}^0_1 \,$. We interpret the relation \eqref{13}$_1$ as a first order linear partial differential equation for the unknown function $ W( {\boldsymbol{E}^e}, {\boldsymbol{K}^e})$, which depends on 12 independent scalar arguments (the 12 components of ${\boldsymbol{E}^e}$ and $ {\boldsymbol{K}^e}$ in the tensor basis $\{\boldsymbol{d}_i^0\otimes\boldsymbol{a}^\alpha\}$). According to the theory of differential equations (see e.g., \cite{Vrabie-2003}, Chap. 6), to solve equation \eqref{13}$_1$ we determine 11 first integrals of the associated system of ordinary differential equations
\begin{equation}\label{15}
    \begin{array}{l}
    \dfrac{\mathrm{d}\boldsymbol{E}^e}{\mathrm{d}s} = \boldsymbol{c}(\boldsymbol{E}^e+\boldsymbol{a}), \qquad \dfrac{\mathrm{d}\boldsymbol{K}^e}{\mathrm{d}s} = \boldsymbol{c}(\boldsymbol{K}^e+\boldsymbol{K}^0).
    \end{array}
\end{equation}
We observe that the functions $\boldsymbol{U}_k$  are first integrals:
\begin{equation}\label{16}
    \begin{array}{l}
   \boldsymbol{U}_1=\boldsymbol{F}^{ T}\boldsymbol{F}=(\boldsymbol{E}^e+\boldsymbol{a})^T(\boldsymbol{E}^e+\boldsymbol{a}),\\
   \boldsymbol{U}_2=\boldsymbol{n}^0\boldsymbol{E}^e= \boldsymbol{d}_3 \boldsymbol{F}, \qquad \boldsymbol{U}_3=\boldsymbol{n}^0\boldsymbol{K}^e.
    \end{array}
\end{equation}
Indeed, we have
\begin{equation*}
    \begin{array}{l}
   \dfrac{\mathrm{d}\boldsymbol{U}_1}{\mathrm{d}s} = \Big(\dfrac{\mathrm{d}\boldsymbol{E}^e}{\mathrm{d}s}\Big)^T(\boldsymbol{E}^e\!+\!\boldsymbol{a}) + (\boldsymbol{E}^{e,T}\!+\!\boldsymbol{a})\dfrac{\mathrm{d}\boldsymbol{E}^e}{\mathrm{d}s} = \boldsymbol{0},\vspace{3pt}\\
    \dfrac{\mathrm{d}\boldsymbol{U}_2}{\mathrm{d}s} = \boldsymbol{n}^0 \dfrac{\mathrm{d}\boldsymbol{E}^e}{\mathrm{d}s}=\boldsymbol{n}^0 [\boldsymbol{c}(\boldsymbol{E}^e\!+\!\boldsymbol{a})] = \boldsymbol{0},
    \end{array}
\end{equation*}
in view of \eqref{15}$_1\,$, and analogously for $\boldsymbol{U}_3$. Finally, another independent first integral is the function
\begin{equation}\label{17}
    \begin{array}{l}
   \boldsymbol{U}_4=\boldsymbol{F}^T\mathrm{Grad}_s\boldsymbol{d}_3 =(\boldsymbol{E}^e\!+\!\boldsymbol{a})^T\boldsymbol{c}\,(\boldsymbol{K}^e\!+\!\boldsymbol{K}^0),
    \end{array}
\end{equation}
since $\frac{\mathrm{d}}{\mathrm{d}s}\, \boldsymbol{U}_4=\boldsymbol{0}$ by virtue of relations \eqref{15}. The functions \eqref{16} and \eqref{17} represent in total 11 scalar independent first integrals. Then, the general solution of the first order  partial differential equation \eqref{13}$_1$ is
$$ W =
      \widetilde W\big(\boldsymbol{F}^{ T}\boldsymbol{F} \, , \,
       \boldsymbol{d}_3 \boldsymbol{F} \, , \, \boldsymbol{n}^0\boldsymbol{K}^e \, , \, \boldsymbol{F}^T\mathrm{Grad}_s\boldsymbol{d}_3 \big),$$
which in view of \eqref{13}$_2$ reduces to \eqref{12}. $\hfill\Box$\smallskip

\textbf{Remark 4.} From  Theorem 3 follows that the strain energy \eqref{12} can be alternatively expressed as a function of the following arguments
\begin{equation}\label{18}
    \begin{array}{l}
    W=\widehat{W}(\boldsymbol{\mathcal{E}},\boldsymbol{\gamma},\boldsymbol{\Psi}), \qquad\boldsymbol{\gamma}= \boldsymbol{d}_3 \boldsymbol{F}= \boldsymbol{n}^0\boldsymbol{E}^e, \\
    \boldsymbol{\mathcal{E}}=\frac{1}{2}\,(\boldsymbol{F}^{ T}\boldsymbol{F} -\boldsymbol{a})= \frac{1}{2}\,\boldsymbol{E}^{e, T}\boldsymbol{E}^e+\mathrm{sym}\, \boldsymbol{E}^e_\parallel\,\,,\\
    \boldsymbol{\Psi}= (\boldsymbol{F}^T\mathrm{Grad}_s\boldsymbol{d}_3 - \mathrm{Grad}_s\boldsymbol{n}^0)-\boldsymbol{\mathcal{E}}\,\mathrm{Grad}_s\boldsymbol{n}^0\\
    \quad\,\, = (\boldsymbol{E}^{e, T}\!\!+\boldsymbol{a})\boldsymbol{c} \boldsymbol{K}^e\!+ [\, \frac{1}{2}\,\boldsymbol{E}^{e, T}\boldsymbol{E}^e\!+ \mathrm{skew} \boldsymbol{E}^e_\parallel\,]\boldsymbol{b},
\end{array}
\end{equation}
where we denote by $\boldsymbol{b}=-\mathrm{Grad}_s\boldsymbol{n}^0$. The tensor $\boldsymbol{\mathcal{E}}$ is a second order symmetric tensor accounting for extensional and in-plane shear strains, $\boldsymbol{\gamma}$ is the vector of transverse shear deformation, and $\boldsymbol{\Psi}$ is a second order tensor for the bending and twist strains.$\hfill\Box$\smallskip

The results \eqref{18} are similar to those presented by Zhilin \cite{Zhilin06} for shells without drilling rotations. The tensors $\boldsymbol{\mathcal{E}}$ and  $\boldsymbol{\gamma}$ coincide with those given in \cite{Zhilin06}, but nevertheless the bending-twist tensor $\boldsymbol{\Psi}$ is different. In this respect,  Zhilin  proposed the tensor
\begin{equation}\label{19}
    \begin{array}{l}
   \boldsymbol{\Phi}
    =\! (\boldsymbol{E}^{e, T}\!\!+\!\boldsymbol{a})\boldsymbol{K}^e_\parallel\!- \![ \frac{1}{2}\,\boldsymbol{E}^{e, T}\!\boldsymbol{E}^e\!+ \mathrm{skew} \boldsymbol{E}^e_\parallel\,]\boldsymbol{c}\,\boldsymbol{b}\\
    \quad\, = [\boldsymbol{F}^T\!(\boldsymbol{d}_3\times\!\mathrm{Grad}_s\boldsymbol{d}_3 ) + \boldsymbol{n}^0 \!\times\!\boldsymbol{b}]+\boldsymbol{\mathcal{E}}(\boldsymbol{n}^0\! \times \! \boldsymbol{b}).
\end{array}
\end{equation}
We consider that the definition of the bending-twist tensor in the form \eqref{18}$_4$ is more appropriate since the relation \eqref{19} introduces an additional (unnecessary) rotation of $\mathrm{Grad}_s\boldsymbol{d}_3$ in the plane $\{\boldsymbol{d}_1,\boldsymbol{d}_2\}$.

From \eqref{18} and \eqref{19} we see that in the linearized theory these deformation tensors reduce to:
\begin{equation*} 
    \begin{array}{l}
    \boldsymbol{\mathcal{E}}\stackrel{.}{=} \mathrm{sym}(\boldsymbol{a}\,\mathrm{Grad}_s\boldsymbol{u}),\quad
    \boldsymbol{\gamma}\stackrel{.}{=}    \boldsymbol{n}^0\mathrm{Grad}_s\boldsymbol{u}+ \boldsymbol{c}\,\boldsymbol{\psi}, \\
    \boldsymbol{\Psi}\stackrel{.}{=}  \boldsymbol{c}\,\boldsymbol{\Phi}
    \stackrel{.}{=}   \boldsymbol{c}\,\mathrm{Grad}_s(\boldsymbol{a}\,\boldsymbol{\psi})+ [\mathrm{skew}(\boldsymbol{a}\,\mathrm{Grad}_s\boldsymbol{u})] \boldsymbol{b},
\end{array}
\end{equation*}
where $\boldsymbol{\psi}$ is the vector of small rotations. One can easily see that $\boldsymbol{\mathcal{E}},\boldsymbol{\gamma}$ and $\boldsymbol{\Psi}$ are independent of the drilling rotation $(\boldsymbol{\psi}\cdot\boldsymbol{n}^0)$. In this case one gets the Reissner-type kinematics of shells \cite{Wisniewski10,Neff_Hong_Reissner08} with 5 degrees of freedom.

\section{Isotropic shells: comparison}\label{sec5}

The isotropic shells without drilling rotations have been investigated in details by Zhilin \cite{Zhilin06}, who determined the form of the strain anergy density $W$ as a quadratic function of its arguments $(\boldsymbol{\mathcal{E}},\boldsymbol{\gamma},\boldsymbol{\Phi})$. Suggested by these results, we consider the following strain energy function for elastic shells without drilling rotations (for the simplified case when the coefficients are independent of $\boldsymbol{K}^0$)
\begin{equation}\label{22}
    \begin{array}{l}
    2\widehat{W}(\boldsymbol{\mathcal{E}},\boldsymbol{\gamma},\boldsymbol{\Psi})=  C[(1\!-\!\nu)\|\boldsymbol{\mathcal{E}}\|^2+ \nu(\mathrm{tr}\,\boldsymbol{\mathcal{E}})^2] 
      +\frac{1}{2}\,C(1\!-\!\nu)\kappa\,\boldsymbol{\gamma}^2 \\ 
      \qquad \qquad \qquad + D[\,\frac{1}{2}\,(1\!-\!\nu)\|\boldsymbol{\Psi}\|^2
      +\frac{1}{2}\,(1\!-\!\nu)\,\mathrm{tr}( \boldsymbol{\Psi}^2)+\nu\,(\mathrm{tr}\,\boldsymbol{\Psi})^2],
    \end{array}
\end{equation}
where $\kappa$ is a shear correction factor. If we insert the expression \eqref{18} of $\boldsymbol{\mathcal{E}},\boldsymbol{\gamma}$ and $\boldsymbol{\Psi}$  into \eqref{22}, then we find the form of $W$ in terms of the strain tensors $(\boldsymbol{E}^e, \boldsymbol{K}^e)$. We observe that the resulting energy density  $W (\boldsymbol{E}^e, \boldsymbol{K}^e)$ is a super-quadratic function of its arguments.
In the case of physically linear shells, when only the quadratic terms in $(\boldsymbol{E}^e, \boldsymbol{K}^e)$ are kept, we obtain the simplified expression of the energy density
\begin{equation}\label{24}
    \begin{array}{l}
    2W( {\boldsymbol{E}^e}, {\boldsymbol{K}^e})= C[\nu(\mathrm{tr}\,\boldsymbol{E}^e_\parallel)^2
    +\,\frac{1-\nu}{2}\,   \mathrm{tr}(\boldsymbol{E}^e_\parallel)^2
    +\,\frac{1-\nu}{2}\,   \mathrm{tr}(\boldsymbol{E}^{e,T}_\parallel\boldsymbol{E}^e_\parallel)]
    +C \,\frac{1-\nu}{2}\,\kappa\,\|\boldsymbol{n}^0\boldsymbol{E}^e\|^2\\
    \qquad \qquad \qquad +D[\mathrm{tr}(\boldsymbol{K}^{e,T}_\parallel\boldsymbol{K}^e_\parallel)- \frac{1-\nu}{2}\,   (\mathrm{tr}\boldsymbol{K}^e_\parallel)^2 - \!\nu \,   \mathrm{tr}(\boldsymbol{K}^e_\parallel)^2].
    \end{array}
\end{equation}

Finally, if we compare the relation \eqref{24} with the general form of the strain energy density for isotropic shells \eqref{6} we find the following values for  $\alpha_k$ and $\beta_k$
\begin{equation}\label{25}
    \begin{array}{l}
    \alpha_1=C\,\nu,\qquad \alpha_2=\alpha_3=C\,\frac{1-\nu}{2}\,,\qquad \alpha_4=C\,\frac{1-\nu}{2}\,\kappa,\\
    \beta_1=D\,\frac{\nu-1}{2}\,,\qquad \beta_2=-D\nu,\qquad \beta_3=D,\qquad \beta_4=0.
    \end{array}
\end{equation}

\textbf{Remark 5.} The coefficients $\alpha_k$ and $\beta_k$ given in \eqref{25} for shells without drilling rotations are different from the values \eqref{8} corresponding to shells with drilling rotations. This indicates that the two types of shells will have different mechanical responses.$\hfill\Box$\smallskip

\textbf{Remark 6.} The conditions \eqref{7} which insure the existence of minimizers are not satisfied by the coefficients \eqref{25} corresponding to shells without drilling rotations since $\alpha_3-\alpha_2=0\,$, $2\beta_1+\beta_2+\beta_3=0$, and $\beta_4=0$. In this case, the strain anergy function \eqref{6} is not uniformly positive definite, and therefore the proof of existence of minimizers is more difficult (in this respect,  see \cite{Neff_plate07_m3as}). This is in accordance to the results presented by Neff \cite{Neff_plate04_cmt,Neff_plate07_m3as} for a plate model derived directly from the 3D equations of Cosserat elasticity. The comparison between the 6-parameter resultant shell theory and the model developed in \cite{Neff_plate04_cmt,Neff_plate07_m3as} has been presented in \cite{Birsan-Neff-JElast-2013,Birsan-Neff-MMS-2013}.$\hfill\Box$

\bigskip\bigskip
\small{\textbf{Acknowledgements.}
The first author (M.B.) is supported by the german state  grant: ``Programm des Bundes und der L\"ander f\"ur bessere Studienbedingungen und mehr Qualit\"at in der Lehre''.

\bibliographystyle{plain}
\bibliography{literatur_Birsan}

\end{document}